\documentclass[12pt]{article}
\usepackage{amsmath,amssymb,latexsym}

\newtheorem{defn}{Definition}[section]
\newtheorem{lemma}[defn]{Lemma}
\newtheorem{ex}[defn]{Example}

\newtheorem{prop}[defn]{Proposition}
\newtheorem{cor}[defn]{Corollary}
\newtheorem{rem}[defn]{Remark}

\newtheorem{conj}{Conjecture}

\newcommand{\h}{{\cal H}}
\newcommand{\mn}{\mathbb N}

\def\range{{\cal R}}

\def\h{{\cal H}}
\def\bp{\noindent{\bf Proof: \ }}
\def\ep{\noindent{$\Box$}}
\def\<{\langle}
\def\>{\rangle}

\newcommand{\seq[1]}{(#1_n)}

\def\newin {\,\kern-0.4em\in\kern-0.15em}
\def\newsubset {\kern-0.2em\subset\kern-0.2em}

\def\v{\vspace{.1in}}
\def\normsn{$\|\!\cdot\!\|$-$SN$}

\newcommand{\neweq}{\overset{\scriptscriptstyle{\nabla}}{=}}

\title{
Can any unconditionally convergent multiplier be transformed to have the symbol $(1)$ and Bessel sequences  by shifting weights?}

{\small 
\author{
 D.\,T. Stoeva$^{a), b)}$\footnote{ std73std@yahoo.com} \ and P. Balazs$^{b)}$\footnote{ peter.balazs@oeaw.ac.at} \\ 
$^{a)}$ Department of Mathematics, \\ University of Architecture, Civil Engineering and Geodesy,\\
Blvd Hristo Smirnenski 1, 1046 Sofia, Bulgaria\\
$^{b)}$ Acoustics Research Institute, \\
Wohllebengasse 12-14, Vienna A-1040, Austria \\
}
}

\begin{document}
\maketitle \pagestyle{myheadings} 
{\footnotetext[1]{This work was supported by  the WWTF project MULAC ('Frame Multipliers: Theory and Application in Acoustics; MA07-025)} }

\parindent0pt
\parskip1ex

\centerline{The paper is dedicated to Hans Feichtinger on the occasion of his 60th birthday.}

\begin{abstract}
Multipliers are operators that combine (frame-like) analysis, a multiplication with a fixed sequence, called the symbol, and synthesis.
The are very interesting mathematical objects that also have a lot of applications for example in acoustical signal processing.
It is known that bounded symbols and Bessel sequences guarantee unconditional convergence.
In this paper we investigate necessary and equivalent conditions for the unconditional convergence of multipliers. 
In particular we show that, under mild conditions,  unconditionally convergent multipliers can be transformed by shifting weights between symbol and sequence, into multipliers with symbol $(1)$ and Bessel sequences. 
\end{abstract}



\section{Introduction}

Multipliers are operators that have the form
\begin{equation} \label{sec:firstmult1} M_{(m_n),(\phi_n),(\psi_n)} f = \sum_{n=1}^\infty m_n\left< f , \psi_n \right> \phi_n,
\end{equation}
where $(\phi_n)$ and   $(\psi_n)$ are sequences in a Hilbert space $\h$ and $(m_n)$ is a scalar sequence, called the symbol.
In \cite{xxlmult1} the known concept of Gabor multipliers \cite{feinow1} was extended to the general frame and Bessel sequences case.

Multipliers are interesting from a mathematical point of view. They have been investigated for Gabor frames \cite{feimulvar1,benepfand07,doetor09}, for fusion frames \cite{Arias2008581}, for generalized frames \cite{rahgenmul10} and $p$-frames in Banach spaces \cite{rahxxlXX}. The concept of multipliers is naturally related to weighted frames \cite{xxljpa1,stoevxxl09} as well as matrix representation of operators \cite{xxlframoper1}. The later is, in particular, important for the numerical solution of operator equations, see e.g. \cite{dafora05,Dahmenetal07}.
Other applications of multipliers are also possible, in particular in acoustics. Multipliers are applied in psychoacoustical modelling \cite{xxllabmask1,labxxl11}, computational auditory scene analysis \cite{wanbro06}, denoising \cite{majxxl10}, sound synthesis \cite{DepKronTor07} or sound morphing \cite{olivtor10}.
For some applications, an approximation of matrices or operators by multipliers is interesting \cite{xxlframehs07,feiham1}.

For Bessel sequences and bounded symbols multipliers are  always well-defined on all of $\h$ with unconditional convergence and bounded \cite{xxlmult1}.
Multipliers can be unconditionally convergent on all of $\h$ for non-Bessel sequences and non-bounded symbols, plenty of examples can be found in \cite{BStable09}.
Multipliers which are well defined on all of $\h$ are always bounded (see Lemma \ref{lemuncb}), but the unconditional convergence is not always guaranteed, see the multiplier $M_{(1),\Phi,\Psi}$ in Example \ref{wdnwd}.
In this paper we  focus on the unconditional convergence of multipliers.

Clearly, the roles of the sequences and the symbols in Equation \eqref{sec:firstmult1} are not independent, some weights can be shifted between those objects. We want to solve the following questions: 
Can we find a {\em `canonical form'} of an unconditional convergent multiplier by shifting weights? 
In particular, as it is known, that a multiplier involving a bounded symbol and Bessel sequences is unconditionally convergent, can we reach such a construction by shifting weights for {\em any} unconditionally convergent multiplier?
Can we connect the invertibility of multipliers to the frame property?
Here we give partial answers and formulate a conjecture for the open question.

In Section \ref{sec:Motiv0} we formulate the questions as motivation for this paper in full details.  
In Section \ref{sec:prel0}, we specify the notation and state the needed results for the main part of the paper. In Section \ref{sec:uncconv}, the unconditional convergence of multipliers is considered; sufficient and equivalent conditions are determined.  
In Section \ref{sec:interpl0}  we give partial answers of the questions posed in Section \ref{sec:Motiv0}. We determine several classes of multipliers, where the Conjecture is true.
Furthermore, we investigate if, by such a shifting, we can also reduce  unconditionally and invertible multipliers to a certain, 'canonical' form. We determine several classes of multipliers which can be reduced to frame multipliers with symbol $(1)$.

\subsection{Motivation} \label{sec:Motiv0} 

In connection to the questions about the re-weighting of symbol and sequence we introduce the following notation:
for sequences $\nu=(\nu_n), \Theta=(\theta_n), \Xi=(\xi_n)$, we will write  
$M_{m,\Phi,\Psi} \neweq  M_{\nu,\Xi,\Theta}$ if there exist scalar sequences $(c_n)$, $(d_n)$ so that $\xi_n=c_n \phi_n$, $\theta_n=d_n\psi_n$, and $m_n=\nu_n c_n \overline{d_n}$ for every $n$.

When $\Phi$ and $\Psi$ are Bessel sequences for $\h$, and  $m\newin\ell^\infty$, then
 $M_{m,\Phi,\Psi}$ is unconditionally convergent on $\h$ \cite{xxlmult1}.

This is only a sufficient condition. For example, 
the multiplier $M_{(n),(\frac{1}{n}e_n),(\frac{1}{n}e_n)}$ is unconditionally convergent on $\h$ and $m=(n)\notin\ell^\infty$.
But note that 
$M_{(n),(\frac{1}{n}e_n),(\frac{1}{n}e_n)}$ can be written as $M_{(1), (e_n), (\frac{1}{n}e_n)}$.
Many examples of unconditionally convergent multipliers $M_{m,\Phi,\Psi}$  with $m\notin\ell^\infty$ or non-Bessel $\Phi$  can be found in \cite{BStable09}. All these multipliers can be transformed into the form $M_{(1),Bessel,Bessel}$ by shifting weights. 

On the other hand, the multiplier $M_{(1),\Phi,\Psi}$ in Example \ref{wdnwd} 
 is well-defined on $\h$, but not unconditionally convergent on $\h$. 
 The sequence $m=(1)$ can not be written in the way $(c_n \overline{d_n})$ so that both $(c_n\phi_n)$ and $(d_n \psi_n)$ are Bessel for $\h$.

The above observations lead to the following question:

\newpage
  \begin{quote} [{$\bf Q_{UC}$}]  {\em If $M_{m,\Phi,\Psi}$ is unconditionally convergent on $\h$,  
 do scalar sequences $\seq[c]$ and $\seq[d]$ exist so that $M_{m,\Phi,\Psi}\neweq M_{(1),(c_n\phi_n),(d_n\psi_n)}$, where 
 $(c_n\phi_n)$  and $(d_n \psi_n)$ are Bessel for $\h$? 
  }
  \end{quote}

The above question is clearly equivalent to the following one:  {\it Does unconditional convergence of $M_{m,\Phi,\Psi}$ on $\h$ imply that $M_{m,\Phi,\Psi}\neweq M_{(\widetilde{m}_n)\in\ell^\infty,(c_n\phi_n),(d_n\psi_n)}$ with  $(c_n\phi_n)$ and $(d_n\psi_n)$ being Bessel for $\h$?}

We can we give a partial answer to $\bf Q_{UC}$: 
\begin{prop} \label{mppuc3} For $M_{m,\Phi,\Psi}$ define the following conditions:
\begin{itemize}
   \item[$\mathcal{P}_1${\rm :} ] $(|m_n|\cdot \|\phi_n\|\cdot \|\psi_n\|)$ is norm-bounded below and $M = M_{m,\Phi,\Psi}$ is unconditionally convergent. 
   \item[$\mathcal{P}_2${\rm :} ] $\exists$ $(c_n)$ and $(d_n)$ so that $M_{m,\Phi,\Psi}\neweq M_{(1),(c_n\phi_n),(d_n\psi_n)}$, where 
 $(c_n\phi_n)$  and $(d_n \psi_n)$ are $\|\!\cdot\!\|$-semi-normalized and Bessel for $\h$.
 \item[$\mathcal{P}_3${\rm :} ] $\exists$ 
 $(c_n)$ and $(d_n)$ so that $M_{m,\Phi,\Psi}\neweq 
 M_{(1),(c_n\phi_n),(d_n\psi_n)}$, 
 where 
 $(c_n\phi_n)$  and $(d_n \psi_n)$ are  Bessel for $\h$.
  \end{itemize}

 For these conditions we have
$\mathcal{P}_1 \Leftrightarrow \mathcal{P}_2 \Rightarrow \mathcal{P}_3$ and $ \mathcal{P}_3 \nRightarrow \mathcal{P}_1$.
\end{prop}

So, under the condition that $(|m_n|\cdot \|\phi_n\|\cdot \|\psi_n\|)$ is norm-bounded below, the question [{$\bf Q_{UC}$}] can be answered affirmatively. 
Also if $§\Phi = \Psi$ we can positively answer the question, see Proposition \ref{q1yes}.

Testing an enormous number of examples of unconditionally convergent multiplier lead us to believe in the following conjecture:
\begin{conj} $M_{m,\Phi,\Psi}$ is unconditionally convergent if and only if $\mathcal{P}_3$ is fulfilled.
\end{conj}
In short, this means that the answer to question [{$\bf Q_{UC}$}] would always be 'Yes'.
\\

By \cite{balsto09new} we know that the invertibility of multipliers is connected to the frame condition of the involved sequences. So in this case we can ask:
\begin{quote} 
[{$\bf Q_{Inv1}$}] {\em
If $M_{m,\Phi,\Psi}$ is unconditionally convergent and invertible,  
 do sequences $\seq[c]$ and $\seq[d]$ exist so that 
$M_{m,\Phi,\Psi}\neweq M_{(1),(c_n\phi_n),(d_n\psi_n)}$, where 
 $(c_n\phi_n)$  and $(d_n \psi_n)$ are frames for $\h$? 
 }
\end{quote}
We determined classes of multipliers for which [$Q_{Inv1}$] has affirmative answer, see Section  \ref{sec:interpl0}.

Note that if the answer of Question [$Q_{UC}$]  is 'Yes', then the answer of [$ Q_{Inv1}$] is also 'Yes', see Section   \ref{sec:interpl0}.
In particular this means that if Conjecture 1 is true, every invertible and unconditionally convergent multiplier can be written as $M_{(1),frame,frame}$.

\section{Notation and preliminaries} \label{sec:prel0}

Throughout the paper 
$\h$ denotes a Hilbert space and $\seq[e]$ denotes an orthonormal basis of $\h$. 
The notion {\it operator} is used for linear mappings. The range of an operator $G$ is denoted by $\range(G)$. 
The identity operator on $\h$ is denoted by $I_\h$. The operator $G:\h\to \h$ is called {\it invertible} if 
 there exists a bounded operator $G^{-1}:\h\to \h$  such that $GG^{-1}=G^{-1}G=I_\h$. 
Throughout the paper, we work with a fixed infinite, but countable index set $J$, and, without loss of generality, $\mn$ is used as an index set, also implicitly.

The notation $\Phi$ (resp. $\Psi$) is used to denote the sequence $\seq[\phi]$ (resp. $\seq[\psi]$) with elements from $\h$; 
 $m$ denotes a 
 complex scalar sequence $\seq[m]$, $\overline{m} = \seq[\overline{m}]$ and $m\Phi = (m_n\phi_n)$.  Recall that $m$ is called {\it semi-normalized} 
  if there exist constants $a,b$ such that $0<a\leq |m_n|\leq b<\infty$, $\forall n$. 
  If $(\|\phi_n\|)$ is semi-normalized, then $\Phi$ is called {\it $\|\!\cdot\!\|$-semi-normalized}. 
    If $\inf_n \|m_n\| >0$ (resp. $\inf_n \|\phi_n\| >0$), the sequence $m$ (resp. $\Phi$) will be called {\it norm-bounded below}, in short $NBB$.

\subsubsection*{Bessel sequences, frames, Riesz bases}
   
Recall that $\Phi$ is called a {\it Bessel sequence} (in short, {\it Bessel}) {\it for $\h$ with bound $B_\Phi$} if $B_\Phi < \infty$ 
and $\sum |\<h,\phi_n\>|^2 \leq B_\Phi\|h\|^2$ for every $h\in\h$. A Bessel sequence $\Phi$ with bound $B_\Phi$  is called a {\it frame for $\h$ with bounds $A_\Phi, B_\Phi$}, if $A_\Phi>0$ and $A_\Phi\|h\|^2\leq  \sum |\<h,\phi_n\>|^2 $ for every $h\in\h$. 
The sequence $\Phi$ is called a Riesz basis for $\h$ with bounds $A_\Phi,B_\Phi$,  
if $\Phi$  is complete in $\h$, $A_\Phi>0$ and 
$A_\Phi \sum |c_n|^2\leq \|\sum c_n \phi_n\|^2\leq B_\Phi \sum |c_n|^2$,
$\forall \seq[c]\in\ell^2$.
Every Riesz basis for $\h$ with bounds $A,B$ is a frame for $\h$ with bounds $A,B$. 
For standard references for frame theory and related topics see \cite{Casaz1,ole1,he98-1}.

\subsubsection*{Unconditional convergence}

  A series $\sum \phi_n$ is called {\it unconditionally convergent} 
 if $\sum \phi_{\sigma(n)}$   converges for every permutation $\sigma(n)$ of $\mn$. 
We will use the following known results about unconditional convergence:

\begin{prop}  \label{hprop} For a sequence $\Phi$, the following statements hold.
\begin{itemize} 
\item[{\rm (i)}] {\rm \cite[Th. 3.16]{he98-1}} If $\sum \phi_n$ converges unconditionally, then
$\sum \|\phi_n\|^2<\infty$. 
\item[{\rm (ii)}] {\rm \cite{he98-1, Orl, Pettis}} The following conditions are equivalent.

$\bullet$ $\sum_n \phi_n$ converges unconditionally.

$\bullet$ Every subseries $\sum_k \phi_{n_k}$ converges.

$\bullet$ Every subseries $\sum_k \phi_{n_k}$ converges weakly.

$\bullet$ $\sum_n \lambda_n \phi_n$ converges for every bounded sequence of scalars $(\lambda_n)$.

\item[{\rm (iii)}] {\rm \cite[Th. 8.3.6]{he98-1}} If $\Phi$ is a Riesz basis for $\h$, then 
$\sum c_n \phi_n$ converges unconditionally  if and only if $\sum c_n \phi_n$ converges.
\item[{\rm (iv)}] If $\Phi$ is a $NBB$ Bessel sequence for $\h$, then $\sum c_n \phi_n$ converges unconditionally if and only if $\seq[c]\in\ell^2$.  

\end{itemize}
\end{prop}
If $\Phi$ is a $NBB$ frame for $\h$, the conclusion of Proposition \ref{hprop}(iv) is proved in \cite[Th. 8.36]{he98-1}. The proof in \cite{he98-1} uses only validity of the upper frame condition, so the property is shown for Bessel sequences.

Concerning Proposition \ref{hprop}(iv), note that if the condition \lq\lq norm-bounded below\rq\rq \ is omitted, then the conclusion does not hold in general, because $\sum c_n \phi_n$ might converge unconditionally for some $\seq[c]\notin\ell^\infty$, see \cite[Ex. 8.35]{he98-1}.

\subsubsection*{Multipliers}
For any $\Phi$, $\Psi$ and any $m$ (called {\it weight} or {\it symbol}), the operator $M_{m,\Phi,\Psi}$, given by $$M_{m,\Phi, \Psi} f=\sum m_n \<f,\psi_n\> \phi_n, \ f\newin\h,$$ is called a {\it multiplier} \cite{xxlmult1}. 
The multiplier $M_{m,\Phi,\Psi}$ is called {\it unconditionally convergent} 
if $\sum m_n\<f,\psi_n\>\phi_n$ converges unconditionally for every $f\in\h$.

Depending on $m,\Phi$, and $\Psi$, the multiplier $M_{m,\Phi,\Psi}$ might not be well defined  (i.e. might not converge for some $f\in\h$) 
 or it might be well defined but not unconditionally convergent. 
First observe that  $M_{m,\Phi,\Psi}$ being well defined on all of $\h$ is not equivalent to $M_{m,\Psi,\Phi}$ 
being well defined on all of $\h$: 

\begin{ex} \label{wdnwd} Let 
$\Phi=(e_1,  e_1,  - e_1,   e_2,   e_1, - e_1,  e_3,  e_1,   - e_1,   \ldots)$ and
 $\Psi=(e_1, e_1, e_1, e_2,   e_2, e_2,  e_3,  e_3,  e_3,   \ldots)$. Then $M_{(1),\Phi,\Psi}=I_\h$  and $M_{(1),\Psi,\Phi}$ is not well-defined. 
\end{ex}

The following statements about well definedness can be easily proved:

 \begin{lemma} \label{lemuncb} For any $\Phi,\Psi$ and $m$, the following holds. 
\begin{itemize}
\item[{\rm (i)}] Let $M_{m,\Phi,\Psi}$ be well defined on all of $\h$. Then $M_{m,\Phi,\Psi}$ is bounded 
 and $M_{\overline{m},\Psi,\Phi}$ equals $M_{m,\Phi,\Psi}^*$ in a weak sense. 
\item[{\rm (ii)}] If $M_{m,\Phi,\Psi}$ and  $M_{\overline{m},\Psi,\Phi}$ are well defined on $\h$, then $M_{\overline{m},\Psi,\Phi}=M_{m,\Phi,\Psi}^*$.
\end{itemize}
\end{lemma}

\section{Necessary and Equivalent Conditions for the Unconditional Convergence of Multipliers}\label{sec:uncconv}

As one can see in Example \ref{wdnwd}, well-definedness of $M_{m,\Phi,\Psi}$ is not equivalent to well-definedness of $M_{m,\Psi,\Phi}$.
If the notion of well-definedness is replaced by the stronger notion of unconditional convergence, then  equivalences hold as follows:

\begin{lemma} \label{uncnew} For any $m$, $\Phi$, and $\Psi$, the following statements are equivalent.
\begin{itemize}
\item[{\rm (i)}] $M_{m,\Phi,\Psi}$ is unconditionally convergent. 
\item[{\rm (ii)}]   $M_{\overline{m},\Psi,\Phi}$ is unconditionally convergent.
\item[{\rm (iii)}]  $M_{m,\Psi,\Phi}$ is unconditionally convergent.
\item[{\rm (iv)}]   $M_{(|m_n|),\Psi,\Phi}$ is unconditionally convergent.

 \end{itemize}
\end{lemma}
\bp (i) $\Leftrightarrow$ (ii): Let $M_{m,\Phi,\Psi}$ be unconditionally convergent.  
By Proposition \ref{hprop}(ii), every subseries $\sum_k m_{n_k} \<f,\psi_{n_k}\> \phi_{n_k}$ converges for every $f\in \h$, which implies that every subseries 
$\sum_k \overline{m}_{n_k} \<g,\phi_{n_k}\> \psi_{n_k}$ converges weakly for every $g\in \h$. Now Proposition \ref{hprop}(ii) implies that $\sum_n \overline{m}_n \<g,\phi_n\> \psi_n$ converges unconditionally for every $g\in\h$.

 (iii) $\Leftrightarrow$ (iv):
Fix $f\in\h$ and assume that $M_{m,\Psi,\Phi}f$ is unconditionally convergent. 
Then every subseries $\sum_k m_{n_k} \<f,\phi_{n_k}\> \psi_{n_k}$ converges unconditionally. 
Consider the sequence $(\lambda_n)$ given by $\lambda_n=\frac{|m_n|}{m_n}$ if $m_n\neq 0$ and $\lambda_n=0$ if $m_n=0$.
Applying Proposition \ref{hprop}(ii) with the bounded sequences $(\lambda_{n_k})_k$, it follows that every subseries $\sum_k |m_{n_k}| \<f,\phi_{n_k}\> \psi_{n_k}$ 
 converges. Now apply again Proposition \ref{hprop}(ii).
 
 The converse follows analogously. 
 
 (ii) $\Leftrightarrow$ (iv) follows from (iii) $\Leftrightarrow$ (iv).
\ep

There exist multipliers which are well defined on all of $\h$ but 
not unconditionally convergent, see $M_{(1),\Phi,\Psi}$ in Example \ref{wdnwd}.
For Bessel sequences and bounded symbols the multiplier is always unconditionally convergent \cite{xxlmult1}. 
Note that this is  only a sufficient condition. Multipliers can be unconditionally convergent 
even in cases when $m\notin\ell^\infty$ or at least one of the sequences is not Bessel. For example, consider $M_{(n^2),(\frac{1}{n}e_n),(\frac{1}{n}e_n)}=I_\h$ and $M_{(1),(\frac{1}{n}e_n),(ne_n)}=I_\h$.
The following statement gives necessary conditions for unconditional convergence:

\begin{prop} \label{lem31} Let $M_{m,\Phi,\Psi}$ be unconditionally convergent. 
\begin{itemize}
\item[{\rm (i)}] 
Then $(m_n\cdot \|\phi_n\| \cdot \psi_n)$ and $(m_n\cdot \|\psi_n\| \cdot \phi_n)$ are Bessel for $\h$.  
\item[{\rm (ii)}] If $\Phi$ ($\Psi$, $m\Phi$, $m\Psi$, respectively) is $NBB$, then $m\Psi$ ($m\Phi$,  $\Psi$, $\Phi$, respectively) is a Bessel sequence for $\h$.
\item[{\rm (iii)}] 
If both $\Phi$ and $\Psi$ are $NBB$, then $m\in\ell^\infty$.
\item[{\rm (iv)}] 
If $\Phi$, $\Psi$ and $m$ are $NBB$, then $m$ is semi-normalized and both $\Phi$ and $\Psi$ are Bessel sequences for $\h$. 
\end{itemize}
\end{prop}
\bp (i)
It follows from Proposition \ref{hprop}(i) that $( \<f,m_n \cdot \|\phi_n\|\cdot \psi_n\> )\in\ell^2$ for every $f\newin\h$. This implies that $(m_n \cdot \|\phi_n\|\cdot\psi_n)$ is Bessel for $\h$.
Now use Proposition \ref{uncnew} and apply what is already proved to $M_{m,\Psi,\Phi}$.

(ii)-(iii) follow easily from (i); (iv) follows from (ii)-(iii).
\ep

\begin{rem} \begin{enumerate}
 \item {\rm Concerning Lemma \ref{lem31}(ii):} if $\Phi$ is not $NBB$, then $m\Psi$ does not need to be a Bessel sequence for $\h$, see \cite[Ex. 4.1.7(i), 4.1.12(i)]{BStable09}. 
 \item {\rm Concerning Lemma \ref{lem31}(iii):} if at least one of $\Phi$ and $\Psi$ is not $NBB$, then $m$ does not need to be in $\ell^\infty$, see \cite[Ex. 4.1.11(i)]{BStable09}.
\end{enumerate}
\end{rem}

For the special cases of Gabor and wavelet systems, which are always $NBB$, Proposition \ref{lem31} gives the following result:
\begin{cor} 
Let $\Phi$ and $\Psi$ be Gabor (or wavelet) systems, and $m$ be $NBB$.
Then $M_{m,\Phi,\Psi}$ is unconditionally convergent if and only if $\Phi$ and $\Psi$ are Bessel for $\h$ and $m$ is semi-normalized.
\end{cor}

Above we have seen sufficient or necessary conditions for the unconditional convergence of multipliers.  Proposition \ref{lem32} and Corollary \ref{lem32a} give conditions which are necessary and sufficient under certain assumptions.

\begin{prop} \label{lem32} For a multiplier $M_{m,\Phi,\Psi}$, the following statements hold.
\begin{itemize}

\item[{\rm (i)}] Let $\Phi$ be a $NBB$ Bessel sequence for $\h$. 
Then 

$M_{m,\Phi,\Psi}$ is unconditionally convergent  $\Leftrightarrow$ 
$m\Psi$ is Bessel for $\h$.

\item[{\rm (ii)}] Let $\Phi$ be a Riesz basis for $\h$. Then

$M_{m,\Phi,\Psi}$ is well defined on $\h$ $\Leftrightarrow$   $M_{m,\Phi,\Psi}$ is unconditionally convergent   
$\Leftrightarrow$  $m\Psi$ is Bessel for $\h$  $\Leftrightarrow$
 $M_{m,\Psi,\Phi}$ is well defined on $\h$
 $\Leftrightarrow$ $M_{m,\Psi,\Phi}$ is unconditionally convergent. 
\item[{\rm (iii)}] Let $\Phi$ be a Riesz basis for $\h$ and $\Psi$ be $NBB$. Then 

$M_{m,\Phi,\Psi}$ (or $M_{m,\Psi,\Phi}$) is well defined on $\h$ 
$\Rightarrow$ $m\in\ell^\infty$. The converse does not hold in general.

\item[{\rm (iv)}] 
If $\Phi$ and $\Psi$ are Riesz bases for $\h$, then 
$M_{m,\Phi,\Psi}$ is well defined on $\h$ if and only if $m\in\ell^\infty$. 
\end{itemize} 
If it is moreover assumed that $m$ is $NBB$ (resp. semi-normalized), then each of the equivalent assertions in (i) and (ii) implies (resp. is equivalent to) $\Psi$ being Bessel for $\h$. 
\end{prop}

\bp 
(i) 
By Proposition \ref{hprop}(iv), $M_{m,\Phi,\Psi}$ is unconditionally convergent  if and only if $(\<f,\overline{m}_n\psi_n\>)\in\ell^2$, $\forall f\in\h$  if and only if
$m\Psi$ is Bessel for $\h$. 

(ii) 
The first equivalence follows from Proposition \ref{hprop}(iii). The second equivalence follows from (i), because Riesz bases are $NBB$ Bessel sequences.

For the third equivalence, consider $M_{m,\Psi,\Phi} f=\sum \<f, \phi_n\> m_n \psi_n,$ $f\newin\h$. 
The sequence $m\Psi$ is Bessel for $\h$ if and only if $\sum c_n m_n\psi_n$ converges for every $\seq[c]\in\ell^2$ if and only if $\sum \<f, \phi_n\> m_n \psi_n$ converges for every $f\newin\h$, because  $\Phi$ is a Riesz basis for $\h$.

To complete the last equivalence, use Proposition \ref{uncnew}.

(iii)  Assume that $M_{m,\Phi,\Psi}$  is well defined, or equivalently, by (ii), that  $M_{m,\Psi,\Phi}$ is well defined.
Let $a_\Psi>0$ denote a lower bound for $(\|\psi_n\|)$.  By (ii),  $m\Psi$ is Bessel for $\h$. Then $a_\Psi |m_n|\leq \|m_n\psi_n\|\leq \sqrt{B_{m\Psi}}$, which implies that $m$ belongs to $\ell^\infty$. 
For the converse, consider the multiplier $M_{(\frac{1}{n}), (e_n),(n^2 e_n)}$, which is not well defined.

(iv)
One of the directions is clear, the  other one follows from (iii).
\ep

\begin{rem} \label{cex1}
{\rm  1. Concerning Prop. \ref{lem32}(i):} If $\Phi$ is Bessel for $\h$, which is non-$NBB$, then the conclusion of Proposition \ref{lem32}(i) might fail. Consider $\Phi=(\frac{1}{2}e_1,e_2,\frac{1}{2^2}e_1,e_3,\frac{1}{2^3}e_1,e_4,\ldots)$, which is Bessel for $\h$, and  $\Psi=(e_1,e_2,e_1,e_3,e_1,e_4,\ldots)$, which is non-Bessel for $\h$. Then $M_{(1),\Phi,\Psi}=M_{(1),\Psi,\Phi}=I_\h$   
with unconditional convergence on $\h$.
 However, $m\Psi=\Psi$ is not Bessel for $\h$.

{\rm 2. Concerning Prop. \ref{lem32}(iii)}:  If $\Phi$ is a Riesz basis for $\h$ and $\Psi$ is non-$NBB$, then well-definedness of  
 $M_{m,\Phi,\Psi}$  does not require $m\in\ell^\infty$.  Consider for example the multiplier $M_{(n), (e_n), (\frac{1}{n}e_n)}$.
 \end{rem}

By Proposition \ref{lem31}(i), a necessary condition for the unconditional convergence of $M_{m,\Phi,\Psi}$ is the sequences $(|m_n|\cdot \|\phi_n\| \cdot \psi_n)$ and $(|m_n|\cdot \|\psi_n\| \cdot \phi_n)$ being Bessel for $\h$. 
It is not difficult to see that this condition is furthermore sufficient under an additional assumption: 

\begin{cor} \label{lem32a}
Let  $(m_n \cdot \|\phi_n\| \cdot \|\psi_n\|)$ be $NBB$. Then the following conditions are equivalent:
\begin{itemize}
\item[{\rm (i)}]
$M_{m,\Phi,\Psi}$ is unconditionally convergent.
\item[{\rm (ii)}]
$(m_n\cdot \|\phi_n\| \cdot \psi_n)$ and $(m_n\cdot \|\psi_n\| \cdot \phi_n)$ are Bessel for $\h$.
\item[{\rm (iii)}] 
 $(|m_n|\cdot \|\phi_n\| \cdot \psi_n)$  and
$(\frac{\phi_n}{\|\phi_n\|})$ are Bessel for $\h$. 
\end{itemize}
\end{cor}

\begin{rem} \label{rem2} 
The $NBB$-property of $(|m_n| \cdot \|\phi_n\| \cdot \|\psi_n\|)$ is not a necessary condition for the unconditional convergence of $M_{m,\Phi,\Psi}$.
 If $(|m_n| \cdot \|\phi_n\| \cdot \|\psi_n\|)$ is non-$NBB$, then unconditional convergence of a multiplier is possible (for example, consider $M_{(\frac{1}{n}),(e_n),(e_n)}$) and non-unconditional convergence of a multiplier is also possible
(for example, consider $\Phi$ and $\Psi$ from Example \ref{wdnwd} and $m=(1,1,1,1,\frac{1}{2}, \frac{1}{2}, 1,1,1,1,\frac{1}{4},\frac{1}{4},\ldots)$, then $M_{m,\Phi,\Psi}=I_\h$  and  $M_{m,\Phi,\Psi}$ is not unconditionally convergent). 
\end{rem}

\subsubsection*{Invertibility and unconditional convergence of multipliers}

 As one can see in Example \ref{wdnwd}, if  $M_{m,\Phi,\Psi}$ is invertible, then  $M_{\overline{m},\Psi,\Phi}$ (resp.  $M_{m,\Psi,\Phi}$) does not need to be neither invertible nor well-defined. But with additional assumptions we can show the following:
 
 \begin{prop} \label{lemuncbinv} For any $\Phi,\Psi$ and $m$, the following holds. 
\begin{itemize}
 \item[{\rm (i)}] Let $M_{m,\Phi,\Psi}$ be invertible  and let $M_{\overline{m},\Psi,\Phi}$ be well defined. Then $M_{\overline{m},\Psi,\Phi}$ is invertible  and $M_{\overline{m},\Psi,\Phi}^{-1}=(M_{m,\Phi,\Psi}^{-1})^*$. 
 \item[{\rm (ii)}]  $M_{m,\Phi,\Psi}$ is unconditionally convergent  and invertible   $\Leftrightarrow$ $M_{\overline{m},\Psi,\Phi}$ is unconditionally convergent and invertible.

\end{itemize}
\end{prop}

\bp (i)  follows from Lemma \ref{lemuncb}(ii).

(ii) follows from Proposition  \ref{uncnew} and Lemma \ref{lemuncb}(ii). 
\ep

As a consequence the following result about dual sequences holds:
\begin{cor} \label{gd} For any $\Phi$ and $\Psi$, the following statements hold.
\begin{itemize} 
\item[{\rm (i)}] If $\sum \<f,\psi_n\> \phi_n =f, \forall f\in\h$, and $\sum \<f,\phi_n\> \psi_n $ converges for every $f\in\h$, then $\sum \<f,\phi_n\> \psi_n =f, \forall f\in\h$.
\item[{\rm (ii)}] $\sum \<f,\psi_n\> \phi_n =f \ \mbox{with unconditional convergence}, \forall f\in\h$, if and only if 
$\sum \<f,\phi_n\> \psi_n =f \ \mbox{with unconditional convergence}, \forall f\in\h$.
\end{itemize}
\end{cor}

\vspace{.1in} Note that Corollary \ref{gd}(ii) generalizes \cite[Lemma 5.6.2]{ole1}, 
which states that if $\Phi$ and $\Psi$ are Bessel sequences, then $\sum \<f,\psi_n\> \phi_n =f, \forall f\in\h$, if and only if $\sum \<f,\phi_n\> \psi_n =f, \forall f\in\h$. In Corollary \ref{gd}(ii) the sequences  $\Phi$ and $\Psi$ do not need to be Bessel sequences for $\h$ - 
for examples with one Bessel and one non-Bessel sequence see \cite[Ex. 4.2.6(i), 4.2.10]{BStable09}, for examples with two non-Bessel sequences see \cite[Ex. 4.1.9(i), 4.1.14(i)]{BStable09}.

\begin{rem} While Lemma \ref{uncnew} gives equivalence of unconditional convergence on $\h$ of $M_{m,\Phi,\Psi}$ and $M_{(|m_n|),\Phi,\Psi}$, note that $M_{m,\Phi,\Psi}$ being unconditionally convergent  and invertible  is not equivalent to $M_{(|m_n|),\Phi,\Psi}$ being unconditionally convergent  and invertible. Consider for example the sequences $\Phi=\Psi=(e_1, e_1, e_2, e_2, e_3, e_3, \ldots)$ and $m=(1, -1, 1, -1, 1, -1, \ldots)$.
\end{rem}

\section{The Interplay of Sequences and Symbols} \label{sec:interpl0}

We have now all necessary tools  
for proving the results in this section.\\

\noindent {\bf Proof of Proposition \ref{mppuc3}:}\\
$\mathcal{P}_1$ $\Rightarrow$ $\mathcal{P}_2$: By $\mathcal{P}_1$ we have that $\phi_n\neq 0$ and $\psi_n\neq 0$, $\forall n\in\mn$.
   By Proposition \ref{lem31}(i), the sequence $(\overline{m_n}\|\phi_n\|\psi_n)$ is Bessel for $\h$. 
Furthermore, $(\overline{m_n}\|\phi_n\|\psi_n)$ is $\|\!\cdot\!\|$-semi-normalized.
  By Corollary \ref{lem32a}, $(\frac{\phi_n}{\|\phi_n\|})$ is Bessel for $\h$. Write $M_{m,\Phi,\Psi}=M_{(1),(\frac{\phi_n}{\|\phi_n\|}),(\overline{m_n}\|\phi_n\|\psi_n)}$. 

 The implications $\mathcal{P}_2$ $\Rightarrow$ $\mathcal{P}_1$ and $\mathcal{P}_2$ $\Rightarrow$ $\mathcal{P}_3$ are clear.
 
 For the implication $ \mathcal{P}_3 \nRightarrow \mathcal{P}_1$, note that $ \mathcal{P}_3$ implies the unconditional convergence, but the $NBB$-property does not necessarily hold, consider for example the multiplier $M_{(1),(\frac{1}{n}e_n),(e_n)}$.
\ep

This means question {$\bf Q_{UC}$} is answered positively 
when $(|m_n|\cdot \|\phi_n\|\cdot \|\psi_n\|)$ is norm-bounded below.
\\

We determine one more class of multipliers, where the answer of {$\bf Q_{UC}$} is affirmative: 
\begin{prop} \label{q1yes} 
The multiplier $M_{m,\Phi,\Phi}$ is unconditionally convergent  if and only if $(\sqrt{m_n}\phi_n)$ is a Bessel sequences for $\h$, where $\sqrt{m_n}$ denotes one (any one) of the two complex square roots of $m_n$, $n\in\mn$.
\end{prop}
\bp
Let $M_{m,\Phi,\Phi}$ be unconditionally convergent. For every $f\in\h$, Lemma \ref{uncnew} 
 implies that 
$M_{(|m_n|),\Phi,\Phi}f$ is unconditionally convergent, which implies that 
$\sum_{n=1}^\infty |m_n| \, |\<f, \phi_n\>|^2 
 <\infty$. 
Assume that $(\sqrt{m_n}\phi_n)$ is not Bessel for $\h$. Then there exists $f\in\h$ so that $(\<f, \sqrt{m_n}\phi_n\>)\notin\ell^\infty$ which contradicts to $\sum_{n=1}^\infty |m_n| \, |\<f, \phi_n\>|^2 
 <\infty$. This completes one of the implications.
 
  The converse implication is clear.
\ep

Hence, for an unconditionally convergent multiplier $M_{m,\Phi,\Phi}$, condition $\mathcal{P}_3$ holds.

\subsubsection*{The Interplay Concerning Invertibility}

 By \cite{balsto09new} we know that the invertibility of multipliers is connected to a frame condition (under some assumptions). Furthermore, we consider the following example:
 
\begin{ex} \label{e1short}
The multiplier $M_{(n),(ne_n),(\frac{1}{n^2}e_n)}$ is  unconditionally convergent  and equal to the Identity operator.
The symbol $m=(n)\notin\ell^\infty$ and the sequences $\Phi$ and $\Psi$ are not frames, but $M_{(n),(ne_n),(\frac{1}{n^2}e_n)}\neweq M_{(1),(e_n),(e_n)}=M_{(1),frame,frame}$.
\end{ex}

On the other hand, observe the following:

\begin{ex}\label{e2} The multiplier
 $M_{(n),(\frac{1}{n}e_n),(\frac{1}{n}e_n)}$ is unconditionally convergent  but not invertible. 
  The sequence $\seq[m]=(n)$ can not be written in the way $(c_n \overline{d_n})$ so that $(c_n\phi_n)$ and $(d_n\psi_n)$ are frames (even, lower frame sequences). Indeed, assume that there exists 
a sequence $(c_n)$ so that $(c_n\phi_n)$ and $(\frac{n}{c_n}\psi_n)$ are 
lower frame sequences with bounds $A_1$ and $A_2$, respectively.
Then
\begin{equation} \label{phif}
A_1\|f\|^2 \leq \sum \mid \<f, c_n \frac{e_n}{n}\>\mid ^2  \mbox{ and } \ A_2\|f\|^2 \leq \sum \mid \<f, \frac{e_n}{c_n} \>\mid ^2, \ \forall f\in\h.\end{equation}
By (\ref{phif}) 
applied with $f=e_j$, $j\in\mn$, it follows that 
$   A_1 j^2\leq \mid c_j\mid ^2\leq \frac{1}{A_2}$, $\forall j\in\mn$, 
 which is a contradiction. 
\end{ex}

Examples \ref{e2} and \ref{e1short} 
lead naturally to the  question [{$\bf Q_{Inv1}$}], which is equivalent to
\begin{quote} 
[{$\bf Q_{Inv\infty}$}] {\em
If $M_{m,\Phi,\Psi}$ is unconditionally convergent  and  invertible, 
 do sequences $(\widetilde{m}_n)\in\ell^\infty$, $\seq[c]$ and $\seq[d]$ exist so that 
  $M_{m,\Phi,\Psi}\neweq M_{(\widetilde{m}_n),(c_n\phi_n),(d_n\psi_n)}$ where 
 $(c_n\phi_n)$  and $(d_n \psi_n)$ are frames for $\h$? 
 }
  \end{quote}

If $M_{m,\Phi,\Psi}$ is invertible, but not unconditionally convergent  (see $M_{(1),\Phi,\Psi}$ in Example \ref{wdnwd}), then $M_{m,\Phi,\Psi}\neweq M_{(\widetilde{m}_n)\in\ell^\infty,Bessel,Bessel}$ is clearly not possible.
\\

\v
Note that if Conjecture 1 is true, then the answer of [{$\bf Q_{Inv1}$}] is always affirmative, as it is connected to [{$\bf Q_{UC}$}].
This is because, by \cite{balsto09new}, if the multiplier $M_{m,\Phi,\Psi}$ is invertible, $m\in\ell^\infty$, and $\Phi$ and $\Psi$ are Bessel for $\h$, then $\Phi$ and $\Psi$ must be frames for $\h$. 
Using this connection we can determine certain classes, as in the unconditional case, where 
we give an affirmative answer of [{$\bf Q_{Inv1}$}]:

\begin{cor} \label{mpp3} Let $M_{m,\Phi,\Psi}$ be  invertible. Define $\mathcal{P}_1$ as in Proposition  \ref{mppuc3} and
\begin{itemize}
      \item[$\widetilde{\mathcal{P}_2}${\rm :} ] $\exists$ $(c_n)$ and $(d_n)$ so that $M_{m,\Phi,\Psi}\neweq M_{(1),(c_n\phi_n),(d_n\psi_n)}$, where 
 $(c_n\phi_n)$  and $(d_n \psi_n)$ are \normsn and frames for $\h$.
 \item[$\widetilde{\mathcal{P}_3}${\rm :} ] $\exists$ 
 $(c_n)$ and $(d_n)$ so that $M_{m,\Phi,\Psi}\neweq M_{(1),(c_n\phi_n),(d_n\psi_n)}$, where 
 $(c_n\phi_n)$  and $(d_n \psi_n)$ are frames for $\h$.
  \end{itemize}
Then the following relations hold:  $\mathcal{P}_1 \Leftrightarrow \widetilde{\mathcal{P}}_2 \Rightarrow \widetilde{\mathcal{P}}_3$ and 
$ \widetilde{\mathcal{P}}_3 \nRightarrow \mathcal{P}_1$.
\end{cor}
\bp  
For the last implication $ \mathcal{P}_3 \nRightarrow \mathcal{P}_1$, 
consider the multiplier  
$M_{(1),\Phi,\Phi}$, where $\Phi=(\frac{1}{\sqrt{2}}e_1,e_2,\frac{1}{\sqrt{2^2}}e_1, e_3,\frac{1}{\sqrt{2^3}}e_1,e_4, \ldots)$.

The rest  follows from Proposition \ref{mppuc3}. 
\ep

\begin{cor} 
If $\Phi$ and $\Psi$ are Gabor (or wavelet) systems, $m$ is $NBB$ and $M_{m,\Phi,\Psi}$ is unconditionally convergent and invertible, then $\widetilde{\mathcal{P}}_3$ holds.
\end{cor}

\begin{cor} Let $M_{m,\Phi,\Psi}$ be unconditionally convergent and invertible. 
If  $\Psi=\Phi$, then $(\sqrt{m_n}\phi_n)$ is a frame for $\h$ (where $\sqrt{m_n}$ denotes one (any one) of the two complex square roots of $m_n$, $n\in\mn$) and thus, $\widetilde{\mathcal{P}}_3$ holds.
\end{cor}

Additionally we can show:
\begin{prop} \label{mpp3b} Let $M_{m,\Phi,\Psi}$ be unconditionally convergent and invertible. 
If  $\Phi$ is minimal, then $\widetilde{\mathcal{P}}_2$ and $\widetilde{\mathcal{P}}_3$ hold.
\end{prop}
\bp 
Let $\Phi$ be minimal. By the invertibility of $M_{m,\Phi,\Psi}$, $\Phi$ is complete in $\h$. Denote by $(\phi_n^b)$ the unique biorthogonal sequence to $\Phi$. 
Then $\phi_n^b=(M_{m,\Phi,\Psi}^{-1})^*(\overline{m}_n\psi_n)$, $\forall n\in\mn$. 
Therefore,
$$1=|\<\phi_n, \phi_n^b\>|\leq \|\phi_n\|  \cdot |m_n| \cdot \|\psi_n\| \cdot \|M_{m,\Phi,\Psi}^{-1}\|, \ \forall n\in\mn.$$
Hence, $(m_n\cdot \|\phi_n\|\cdot \|\psi_n\|)$ is $NBB$.
Now the unconditional convergence of $M_{m,\Phi,\Psi}$ and Corollary \ref{mpp3} complete the proof.
\ep

 {\bf Acknowledgments} The authors are thankful to H. Feichtinger and D. Bayer for their valuable comments. 
The first author is grateful for the hospitality of the Acoustics Research Institute and the support from the MULAC-project. 
She is also grateful to the University of Architecture, Civil Engineering and Geodesy, and the Department
of Mathematics of UACEG for supporting the present research.

{\small
\bibliographystyle{plain}

 \end{document}